\documentclass[english,10pt]{amsart}

\usepackage{amssymb}
\usepackage[all]{xy}

\usepackage{graphicx}

\usepackage{amsmath,amssymb,enumerate,bbm}

\usepackage[T1]{fontenc}
\usepackage[utf8]{inputenc}
\usepackage[all]{xy}

\usepackage{babel}
\usepackage{amstext}
\usepackage{amsmath}
\usepackage{amsfonts}
\usepackage{latexsym}
\usepackage{ifthen}

\usepackage{xypic}
\xyoption{all}
\newcommand{\rk}{{\rm rk} }

\newtheorem{lemma1}{}[section]

\newenvironment{lemma}{\begin{lemma1}{\bf Lemma.}}{\end{lemma1}}
\newenvironment{example}{\begin{lemma1}{\bf Example.}\rm}{\end{lemma1}}
\newenvironment{examples}{\begin{lemma1}{\bf Examples.}\rm}{\end{lemma1}}

\newenvironment{theorem}{\begin{lemma1}{\bf Theorem.}}{\end{lemma1}}
\newenvironment{proposition}{\begin{lemma1}{\bf Proposition.}}{\end{lemma1}}

\newenvironment{remark}{\begin{lemma1}{\bf Remark.}\rm}{\end{lemma1}}
\newenvironment{remarks}{\begin{lemma1}{\bf Remarks.}\rm}{\end{lemma1}}
\newenvironment{definition}{\begin{lemma1}{\bf Definition.}}{\end{lemma1}}

\newenvironment{conjecture}{\begin {lemma1}{\bf Conjecture.}}{\end{lemma1}}
\newenvironment{question}{\begin {lemma1}{\bf Question.}}{\end{lemma1}}
\newenvironment{assumption}{\begin{lemma1}{\bf Assumption.}}{\end{lemma1}}

\newenvironment{remark*}{{\bf Remark.}}{}
\newenvironment{remarks*}{{\bf Remarks.}}{}
\newenvironment{example*}{{\bf Example.}}{}

\newcommand{\C}{\ensuremath{\mathbb{C}}}
\newcommand{\N}{\ensuremath{\mathbb{N}}}
\newcommand{\PP}{\ensuremath{\mathbb{P}}}

\newcommand{\Homsheaf} { \ensuremath{ \mathcal{H} \! om}}

\newcommand{\merom}[3]{\ensuremath{#1:#2 \dashrightarrow #3}}

\newcommand{\holom}[3]{\ensuremath{#1:#2  \rightarrow #3}}
\newcommand{\fibre}[2]{\ensuremath{#1^{-1} (#2)}}

\makeatletter
\ifnum\@ptsize=0  \fi \ifnum\@ptsize=2  \fi 

\newcommand\sF{{\mathcal F}}
\newcommand\sG{{\mathcal G}}

\newcommand\sI{{\mathcal I}}

\newcommand\sL{{\mathcal L}}
\newcommand\sO{{\mathcal O}}

\DeclareMathOperator*{\sing}{sing}
\DeclareMathOperator*{\reg}{reg}

\newcommand{\chow}[1]{\ensuremath{\mathcal{C}(#1)}}

\newcommand{\Chow}[1]{\ensuremath{\mathcal{C}(#1)}}

\newcommand{\Hilb}{\ensuremath{\mathcal{H}}}
\newcommand{\droite}{\ensuremath{\mbox{line}} }

\title{Algebraic foliations defined by quasi-lines} 
\date{July 27, 2009}

\author{Laurent Bonavero}
\author{Andreas H\"oring}

\address{
Laurent Bonavero,
Institut Fourier, UMR 5582 du CNRS, 
Universit{\'e} de Grenoble I,
BP~74, 38402 Saint-Martin d'H{\`e}res, France 
}
\email{laurent.bonavero@ujf-grenoble.fr}

\address{
Andreas H\"oring, Universit\'e Paris 6, Institut de Math\'ematiques de Jussieu, Equipe de Topologie et G\'eom\'etrie Alg\'ebrique, 175, rue du Chevaleret, 75013 Paris, France
}
\email{hoering@math.jussieu.fr}


\begin{document}

\begin{abstract}
Let $X$ be a projective manifold containing a quasi-line $l$.
An important difference between quasi-lines and lines in the projective space is that 
in general there is more than one quasi-line passing through two given general points.
In this paper we use this feature to construct an algebraic foliation associated to 
a family of quasi-lines. We prove that if the singular locus of this foliation is not 
too large, it induces
a rational fibration on $X$ that maps the general leaf
of the foliation onto a quasi-line in a rational variety.
\end{abstract}

\maketitle


\section{Introduction}


\subsection{Motivation}

Let $X$ be a complex quasiprojective manifold of dimension $n$.
A quasi-line $l$ in $X$ is a {\it smooth} rational curve $f~: \PP^1 
\hookrightarrow X$
such that $f^*T_X$ is the same as for a line in $\PP^n$, {\em i.e.}\ 
is isomorphic to 
$$
{\mathcal O}_{\PP^1}(2) \oplus {\mathcal O}_{\PP^1}(1)^{\oplus n-1}.
$$
Quasi-lines have some of the deformation properties of lines, but there are important differences:
for example if $x$ and $y$ are general points in $X$ there exist only finitely many deformations of $l$ passing through the two points,
but in general we do not have uniqueness\footnote{We denote by $e(X,l)$ the number of quasi-lines through two general points,
see Definition \ref{definitioneinvariant} for a formal definition.}. 
It is now well established that given a variety
$X$ with a quasi-line $l$, the deformations and degenerations
of $l$ contain interesting information on the global geometry of $X$. 
Here is an example of such a result, due to Ionescu and Voica.

\begin{theorem} \cite[Thm.1.12]{IV03} \label{theoremvoica} Let $X$ be a projective manifold containing
a quasi-line $l$. Assume there exists a divisor $D$ such that
$D \cdot l=1$ and $h^0(X,\mathcal O _X(D))=s+1 \geq 2$.
Then there exists a small deformation $l'$ of $l$,
a finite composition of smooth blow-ups $\sigma : \tilde X \to X$  
with smooth centers disjoint from $l'$ and a surjective
fibration $\varphi : \tilde X \to \mathbb P ^s$ with rationally connected general fibre such that
$\varphi$ maps isomorphically $\sigma^{-1}(l')$ to
a line in $\mathbb P ^s$.\footnote{In order to simplify the statements, we'll
simply say that there exists a rational fibration
\merom{\varphi}{(X,l)}{(\PP^s,\droite)}.}  
\end{theorem} 

A disadvantage of this statement is that {\em a priori} there seems to be no relation between the geometry
of the quasi-line $l$ and the existence of the divisor $D$. The goal of this paper is to fill this gap by a construction inspired
by the theory of complex projective manifolds $X$ swept out by linear spaces:
these have been studied for more than twenty years (see \cite{Ein85, ABW92, Sat97, NO07}) and 
an observation common to all these papers is that if the codimension of the linear space is small, then either $X$ is special
(a projective space, hyperquadric {\em etc.}) or it admits a 
fibration such that the fibres are linear spaces.
A powerful tool in their theory is the family of 
lines contained in the linear spaces.
The guiding philosophy of this paper is that the rich geometry of a family of quasi-lines can be used to construct
a natural family of subvarieties that induces a (rational) fibration on $X$.

\subsection{Setup and main results}

Let $X$ be a projective manifold of dimension $n$ containing
a quasi-line $l$.
The main tool used in this paper is an intrinsic foliation $\sF_x$ 
associated to the quasi-lines passing through a general point $x$
of $X$. 
In case the foliation has rank $n-1$, its 
leaves are natural candidates
to play the role of the divisor $D$ in Theorem \ref{theoremvoica}. 
The foliation $\sF_x \subset T_X$ is defined by the following heuristic principle: 
\begin{quotation}
``for $y$ general in $X$, the (closure of the) $\sF_x$-leaf
through $y$ is the smallest subvariety $V \subset X$ containing $y$ and
such that for every $z$ in $V$, every quasi-line through $x$ and $z$
is entirely contained in $V$''.
\end{quotation}  
In a more technical language (see Section \ref{sectionfoliation}) we prove the following theorem.

\begin{theorem} \label{theoremexistencefoliation}
Let $X$ be a projective manifold containing
a quasi-line $l$ and let
$\Hilb_x \subset \chow{X}$ be the scheme 
parametrising deformations and degenerations
of $l$ passing through a general point
$x \in X$. 
Then there exists a unique saturated algebraic 
foliation $\sF_x \subset T_X$ such that
for every general point $y \in X$, the unique $\sF_x$-leaf (cf. Defn. \ref{definitionfoliation})
through $y$ is the minimal $\Hilb_x$-stable projective subvariety through $y$.
\end{theorem}

If $l$ is a line or more generally if 
$e(X,l)=1$, the foliation $\sF_x$ has rank one: 
the leaf through a general point $y$ is the unique quasi-line passing through $x$ and $y$.
This leads immediately to the following question.

\begin{question} \label{questionmain}
Let $X$ be a projective manifold containing
a quasi-line $l$. Let $x$ be a general point in $X$, and 
denote by $\sF_x$ the corresponding foliation.
Can we construct a rational fibration 
$$
\merom{\varphi}{(X,l)}{(Y,l':=\varphi(l))}
$$ 
onto a projective variety $Y$ such that 
\begin{itemize}
\item $l'$ is a quasi-line with $e(Y,l')=1$, and
\item the general $\sF_x$-leaves are preimages
of deformations of $l'$?
\end{itemize}
\end{question}

Suppose for a moment that such a fibration exists: fix two general points 
$x$ and $y$ in $X$, and denote by $\sF_x$ and $\sF_y$ the corresponding foliations.
By hypothesis the unique  $\sF_x$-leaf through $y$ is the preimage of
a quasi-line  through $\varphi(x)$ and $\varphi(y)$.
Analogously the unique  $\sF_y$-leaf through $x$ is the preimage 
of a quasi-line through $\varphi(y)$ and $\varphi(x)$.
Since both quasi-lines are deformations of $l'$  passing through two given general points
the condition $e(Y, l')=1$ implies that they are identical.
Hence the two leaves are identical. More formally we have a natural necessary condition 
for the existence of the fibration.

\begin{assumption}
\label{assumption}
Let $X$ be a projective manifold containing
a quasi-line $l$.
Let $x$ and $y$ be two general points in $X$, and 
denote by $\sF_x$ and $\sF_y$ the corresponding foliations.
Denote by $F_{x,y}$ the unique  $\sF_x$-leaf through $y$
and by $F_{y,x}$ the unique $\sF_y$-leaf through $x$.
Then  we have 
$$
F_{x,y}=F_{y,x}.
$$ 
\end{assumption}

Since every $\sF_x$-leaf contains $x$, the singular locus $\sF^{\sing}_x$ of the foliation 
contains $x$ and it is rather optimistic to expect general leaves to be smooth. 
Our first observation is that under the Assumption \ref{assumption} we have some control on the singularities around $x$.

\begin{proposition}
Under the Assumption \ref{assumption}, 
let $\sL$ be a general $\sF_x$-leaf
and $l \subset X$ be a general quasi-line through $x$ such that
$l \subset \sL$.
Then $l$ is contained in the smooth locus of $\sL$ and is a quasi-line in $\sL$.
In particular $\sL$ is smooth in $x$.
\end{proposition}

Now that we have some local information about the general leaf $\sL$, we can try to understand the global geometry of $\sL$ and $X$.
On the one hand we observe that
$$
\dim_x \sF_x^{\sing} \geq \rk \sF_x-1,
$$
so this looks like a rather tough task.
On the other hand we know that the singular locus of 
the foliation given by lines in the projective space has dimension zero. Thus
if we had a rational fibration $\merom{\varphi}{(X,l)}{(\PP^{n-\rk \sF_x+1}, \droite)}$,
the singular locus would be of dimension exactly $\rk \sF_x-1$.
Our main theorem shows that these necessary conditions
are sufficient to give an affirmative answer
to Question \ref{questionmain}.

\begin{theorem}  \label{maincaseI}  
Let $X$ be a projective manifold containing
a quasi-line $l$. Let $x$ be a general point in $X$, and 
denote by $\sF_x$ the corresponding foliation with general leaf $\sL$.
Assume that $\rk \sF_x<\dim X$ and the Assumption \ref{assumption} holds.
Then the following holds.
\begin{enumerate}
\item[a)] If the foliation $\sF_x$ satisfies
\[
\dim_x \sF_x^{\sing} < \rk \sF_x-1+ \frac{1}{2} (n-\rk \sF_x),
\]
there exists a rational fibration \merom{\varphi_\sL}{(\sL, l)}{\PP^1} such that $l$ is a section\footnote{More precisely,
there exists a finite composition of 
smooth blow-ups $\sigma : \tilde{\sL} \to \sL$  
with smooth centers disjoint from a deformation $l'$
of $l$  and a surjective
fibration $\varphi_\sL : \tilde{\sL} \to \mathbb P ^1$ 
such that $\sigma^{-1}(l')$ is a section of $\varphi_\sL$.}.
\item[b)] If the foliation $\sF_x$ satisfies
\[
\dim_x \sF_x^{\sing} = \rk \sF_x-1,
\]
there exists a rational fibration 
$\merom{\varphi}{(X,l)}{(Y,l')}$ onto 
a projective variety $Y$ of dimension $n-\rk \sF_x+1$
such that $l'$ is a quasi-line with $e(Y, l')=1$.
Moreover the restriction of $\varphi$ to $\sL$ is the fibration $\varphi_\sL$.
\end{enumerate}
\end{theorem} 

Since the quasi-line $l' \subset Y$ satisfies $e(Y, l')=1$, the variety $Y$ is rational by \cite[Prop.3.1]{IN03}. 
In view of Theorem \ref{theoremvoica} it would be nice to know when the quasi-line $l' \subset Y$ identifies 
to a line in a projective space. One guess is the following.

\begin{conjecture}
Let $X$ be a projective manifold containing
a quasi-line $l$ with $e(X,l)=1$. If the (rank one) foliation  
$\sF_x$ satisfies
\[
\dim_x \sF_x^{\sing} = 0,
\]
there exists a birational map $X \dashrightarrow \PP^n$ that is an isomorphism in a neighbourhood
of $l$ and maps $l$ onto a line.
\end{conjecture}

A similar statement 
is shown in \cite[Thm.1.5]{IR07}, but they make the additional hypothesis that all the deformations of $l$ are smooth at $x$.
If the rank of $\sF_x$ equals  $n-1$ all these consideration boil down to a very simple statement.

\begin{proposition} \label{propositioncorankone}
Let $X$ be a projective manifold containing
a quasi-line $l$.
Let $x$ be a general point in $X$, and 
denote by $\sF_x$ the corresponding foliation.
Assume that $\rk \sF_x=n-1$ and the Assumption \ref{assumption} holds.
Then the following holds.
\begin{enumerate}
\item[a)] There exists a rational fibration $\merom{\varphi}{(X,l)}{(\PP^2, \droite)}$.
\item[b)] If the Picard number of $X$ is one, we have $X \simeq \PP^2$. 
\end{enumerate}
\end{proposition}

Although we do not have an example where the Assumption \ref{assumption} fails, it is of course
legitimate to ask what happens in this situation. 
We don't have much hope to understand the geometry of $\sF_x$ or $X$,
but we can address another 
basic problem in the study of quasi-lines: 
computing, or at least bounding the invariant $e(X,l)$. 

\begin{theorem} \label{maincaseII}
Let $X$ be a projective manifold of dimension three 
containing a quasi-line $l$.
Let $x$ be a general point in $X$, and 
denote by $\sF_x$ the corresponding foliation.
Assume that $\rk \sF_x=2$ and the Assumption \ref{assumption} {\bf fails}.
Then 
\[
e(X,l) \leq 16 \frac{(\deg l)^3}{\deg X},
\]
where all the degrees are taken with respect to an 
ample line bundle $H$ on $X$.
\end{theorem}

If the Assumption \ref{assumption} holds it is much more difficult to find a reasonable bound for the invariant $e(X,l)$, 
cf. the discussion after Proposition \ref{propositionbadestimate}.

{\bf Acknowledgements.} The second named author wants to thank Kristina Frantzen for explaining to him 
the nice combinatorial argument in Example \ref{examplepardinicount}.

\subsection{Notation and basic definitions}

We work over the complex field $\C$, topological 
notions always refer to the Zariski topology.
A variety is an integral scheme of finite type over $\C$, 
a manifold is a smooth variety.
A fibration is a surjective 
morphism \holom{\varphi}{X}{Y} between normal varieties such that
$\dim X>\dim Y>0$ and $\varphi_* \sO_X \simeq \sO_Y$,
that is all the fibres are connected. Fibres are 
always scheme-theoretic fibres.
For general definitions we refer to Hartshorne's book \cite{Har77}, 
we will also use the standard terminology
of Mori theory and deformation theory as explained 
in \cite{Deb01, Kol96}.

Let $X$ be a projective variety and let 
$V \subset X$ be a projective subvariety.
Identify 
$V$ to its fundamental 
cycle\footnote{Throughout the whole paper, 
we will not distinguish between an effective cycle and its support.}.
We denote by $[V]$ the point in 
the Chow scheme $\chow{X}$ corresponding to $V$. 

\begin{definition}\label{definitiongeneral}
Let $N \subset \chow{X}$ be a finite union of subvarieties 
parametrising $d$-dimensional cycles $V \subset X$. 
We say that a property holds for a general ({\em resp.} very general) cycle $V$ if 
there exist an {\em open dense} ({\em resp.} dense) subset $N^0 \subset N$ 
such that the property holds for every cycle parametrised by a 
point $[V] \in N^0$.
\end{definition}

For the convenience of the reader we recall some well known fact from  deformation theory:
let $X$ be a complex projective manifold of dimension $n$, 
and let $l \subset X$ be a quasi-line. Since 
$$
N_{l/X} \simeq \sO_{\PP^1}(1)^{\oplus n-1}, 
$$
the  Chow scheme $\chow{X}$ is smooth of dimension $2n-2$ at $[l]$. 
Therefore there exists
a unique irreducible component $\Hilb$ of $\chow{X}$ 
containing $[l]$. Furthermore a curve of $X$ 
corresponding to a general point of $\Hilb$ is a quasi-line.
Let $\Gamma \subset \Hilb \times X$ be the universal family, then
the natural map
$\Gamma \rightarrow X$ is surjective.
Thus for a general point
$x \in X$, the subscheme 
$\Hilb_x \subset \Hilb$ parametrising deformations and degenerations
of $l$ passing through $x$ has pure dimension $n-1$. 
Furthermore the points $[l']$ corresponding 
to quasi-lines are dense in $\Hilb_x$.
Since for any such quasi-line $l' \subset X$  
$$
N_{l'/X} \otimes \sI_x \simeq \sO_{\PP^1}^{\oplus n-1}, 
$$
we see that all the irreducible components 
of $\Hilb_x$ are generically smooth. Setting $\Gamma_x$ for the universal family
we get the following basic diagram:
\begin{equation}
\label{basicdiagram}
 \xymatrix{ 
{\Gamma}_x \ar[d]^{q}\ar[r]^{p} & X \\
\Hilb_x & } 
\end{equation}

\begin{remarks}
\begin{enumerate}
\item[a)] We fix the family $\Hilb$ for the rest of the paper. 
Thus if we say that a property holds for every quasi-line 
through $x$ and $y$, we mean 
every quasi-line through $x$ and $y$ parametrised by $\Hilb$, 
{\em i.e.}\  being a member of the family we are interested in. 
\item[b)] Let us also point out what the meaning of 
general (see Definition \ref{definitiongeneral}) means in this context: 
a property holds for a general quasi-line through $x$ if 
there exists an {\em open dense} subset $\Hilb^0_x \subset \Hilb_x$ such that
it is satisfied by every quasi-line parametrised by $\Hilb^0_x$.
\end{enumerate}
\end{remarks}

We can now give the technical definition of 
the invariant $e(X,l)$. This number has been introduced
by Ionescu and Voica and gives the number of quasi-lines which are deformations
of $l$ and pass through two given general points of $X$. 

\begin{definition}
\label{definitioneinvariant}
Let $X$ be a complex projective manifold of dimension $n$ 
and let $l \subset X$ be a quasi-line.
Let $x \in X$ be a general point and let \holom{p}{{\Gamma}_x}{X}
be the morphism from Diagram \eqref{basicdiagram}. Then we define
\[
e(X,l):= \deg(p).
\]
\end{definition}

Let us finally define the kind of foliations we'll be interested in.

\begin{definition} \label{definitionfoliation}
Let $X$ be a projective manifold. 
A reflexive subsheaf $\sF \subset T_X$ 
is called a foliation on $X$ if it satisfies the following
two conditions.
\begin{enumerate}
\item[a)] The subsheaf $\sF$ is saturated, 
that is
the quotient $T_X/\sF$ is torsion-free. 
Since a torsion-free sheaf is locally free in codimension one, 
this implies that
the regular locus of the foliation, {\em i.e.}\  
the open set
\[
\sF^{\reg} = \{ y \in X \ | \ \sF \subset T_X  \ \mbox{is a subbundle 
in an analytic neighbourhood of $y$} \}
\]
is the complement of a subset of codimension at least two. 
Equivalently the 
singular locus of the foliation
\[
\sF^{\sing} := X \setminus \sF^{\reg}
\]
has codimension at least two.
\item[b)] For every $y \in \sF^{\reg}$, there exists a 
projective subvariety $F_y$ passing through $y$ such that
$$
\sF|_{ \sF^{\reg} \cap F_y} = T_{F_y}|_{ \sF^{\reg} \cap F_y} \subset T_{X}|_{ \sF^{\reg} \cap F_y}.
$$
We call $F_y$ the unique $\sF$-leaf through $y$. 
\end{enumerate} 
\end{definition}

\begin{remarks}
\begin{enumerate}
\item[a)] Our definition of foliation is much stronger than the usual 
one in differential geometry since we want the leaves to be 
proper subvarieties of $X$. Such foliations are also called
{\em algebraic foliation}. 
\item[b)] We also call {\em leaf} of a foliation the closure of the usual leaf
in differential geometry, hoping there will be no confusion for the reader.   
The second condition then implies that the 
singular locus of a leaf $F_y$ is contained in the 
singular locus of the foliation.
\end{enumerate}
\end{remarks}

\section{The foliations associated to a family of quasi-lines} \label{sectionfoliation}

\subsection{Construction of the foliations}

Let $X$ be a complex projective manifold of dimension $n$, 
and let $l \subset X$ be a quasi-line.
Let $x \in X$ be a general point and let
\[
 \xymatrix{ 
{\Gamma}_x \ar[d]^{q}\ar[r]^{p} & X \\
\Hilb_x & } 
\]
be the basic diagram \eqref{basicdiagram}.

Let now $\tilde{\Hilb}_x \subset \Hilb_x$ 
be the maximal smooth open dense subset 
such that the points of $\tilde{\Hilb}_x$ parametrise 
quasi-lines through $x$.
We denote by $\tilde{\Gamma}_x$ the universal family 
in $\tilde{\Hilb}_x \times X$.  
Note that $\tilde{\Gamma}_x$ is a quasi-projective,
not necessarily connected, smooth scheme
that is a $\PP^1$-bundle over $\tilde{\Hilb}_x$.

Let $\tilde{\Gamma}^*_x \subset \tilde{\Gamma}_x$ be the maximal open 
dense subset
such that the restriction of \holom{p}{\Gamma_x}{X} 
to $\tilde{\Gamma}^*_x$ is {\em flat and finite}. 
Up to replacing $\tilde{\Hilb}_x$ by $q(\tilde{\Gamma}^*_x)$, 
we can suppose that the restriction of $q$ to $\tilde{\Gamma}^*_x$
is equidimensional of relative dimension $1$ and surjective.

Denoting by $\tilde{p}$ and $\tilde{q}$  
the restrictions of the natural maps $p$ and $q$ to 
$\tilde{\Gamma}^*_x$, we get a new basic diagram:
\begin{equation}
\label{basicdiagramstar}
 \xymatrix{ 
{\tilde{\Gamma}^*_x} \ar[d]^{\tilde{q}}\ar[r]^{\tilde{p}} & X \\
{\tilde{\Hilb}_x} & } 
\end{equation}

\begin{remarks}
\begin{enumerate}
\item[a)] Since all the quasi-lines parametrised 
by $\Hilb_x$ pass through $x$, 
the map $p: \Gamma_x \rightarrow X$ contracts a divisor onto $x$. 
Therefore the point $x$ is not in the image 
of $\tilde{p}: \tilde{\Gamma}^*_x \rightarrow X$.

\item[b)] If $y$ is a general point in $X$, 
all the quasi-lines through $x$ and $y$ parametrised by $\Hilb$ 
are parametrised by a point of $\tilde{\Hilb}_x$: if this was 
not the case, we could vary $y$ in a $n$-dimensional quasi-projective family
and get an $n-1$-dimensional subscheme of $\Hilb_x$ 
that is disjoint from $\tilde{\Hilb}_x$. 
Since $\tilde{\Hilb}_x$ is dense in $\Hilb_x$ and all 
the irreducible components have dimension $n-1$, we would get a contradiction.
In particular $\tilde{p}$ is dominant and we have
\[
\deg \tilde{p} = \deg p = e(X,l).
\]
\end{enumerate}
\end{remarks}

\begin{example}
Set $X:=\PP^2$, and consider the family of lines on $\PP^2$. 
For $x \in X$ arbitrary, we have $\Hilb_x \simeq \PP^1$ and 
the universal family $\Gamma_x$ 
is the first Hirzebruch surface $\mathbb{F}_1$. All the curves 
parametrised by $\Hilb_x \simeq \PP^1$ are lines through $x$, so
$\Hilb_x=\tilde{\Hilb}_x$ and $\Gamma_x=\tilde{\Gamma}_x$.
The natural map $p: \mathbb{F}_1 \rightarrow \PP^2$ is given 
by the contraction of the unique $(-1)$-curve $E$, so 
$\tilde{\Gamma}^*_x = \mathbb{F}_1 \setminus E$.
\end{example}

\begin{definition}
Let $V \subset X$ be a projective subvariety such that $V$
intersects the image of $\tilde{p}: \tilde{\Gamma}^*_x \rightarrow X$.
We say that $V$ is $\Hilb_x$-stable if
\[
\overline{\tilde{p}(\fibre{\tilde{q}}{\tilde{q}(\fibre{\tilde{p}}{V})}} = V.
\]
Let $y \in X$ be a general point. A projective 
subvariety of $X$ containing $y$ is a minimal
$\Hilb_x$-stable subvariety through $y$ if it is contained in
any $\Hilb_x$-stable subvariety containing $y$.
\end{definition}

\begin{remark}
It is immediate from the definition that for 
any general point $y \in X$, a minimal 
$\Hilb_x$-stable subvariety through $y$ is unique. 
Furthermore a quasi-line passing through $x$ and a 
general point of a $\Hilb_x$-stable subvariety $V$ is contained in $V$.
\end{remark}

We are now ready to state and prove our first result.

{\bf Theorem \ref{theoremexistencefoliation}.}
{\em Let $X$ be a projective manifold of dimension $n$ containing
a quasi-line $l$, and let
$\Hilb_x \subset \chow{X}$ be the scheme 
parametrising deformations and degenerations
of $l$ passing through a general point
$x \in X$. 
Then there exists a unique saturated algebraic 
foliation $\sF_x \subset T_X$ such that
for every general point $y \in X$, the unique $\sF_x$-leaf
through $y$ is the minimal $\Hilb_x$-stable projective
subvariety through $y$.
}

The construction of $\sF_x$ is roughly the same as the one given 
by Ein-K\"uchle-Lazarsfeld \cite{EKL95} or Hwang-Keum \cite{HK03}. 
The basic idea is very simple: let $y \in X$ be a general point
and choose a quasi-line $l$ from $x$ to $y$. 
Let $Z_l$ be the set of quasi-lines through $x$ meeting a general 
point of $l$, {\em i.e.}\  a point
in $l \cap {\rm Im }(\tilde{p})$.
Then there are two possibilities: the curves 
parametrised by $Z_l$ dominate a surface $l \subset S \subset X$  
or map into $l$.
In the second case we have finished, the curve $l$ is 
the minimal $\Hilb_x$-stable subvariety through $y$.
In the first case we choose an irreducible component of the surface 
$S$ and restart the construction, {\em i.e.}\  we 
consider the set $Z_S$ of quasi-lines through $x$ 
meeting a general point of $S$, {\em etc}.

\begin{proof}
Let $x$ be a general point of $X$.

{\it Step 1: construction of the minimal $\Hilb_x$-stable 
subvarieties.} Fix $y \in X$ a general point. 
We define a sequence of subvarieties as follows: 
set $V_0:=y$, and for $i \in \{ 1, \ldots, n+1 \}$,
let $V_i$ be an irreducible component of
\[
\overline{\tilde{p}(\fibre{\tilde{q}}{\tilde{q}(\fibre{\tilde{p}}{V_{i-1}})}}
\]
that has maximal dimension. By flatness 
of $\tilde{p}$, every irreducible component of $\fibre{\tilde{p}}{V_{i-1}}$
has dimension $\dim V_{i-1}$ and dominates 
$V_{i-1}$. Thus we have
\[
V_{i-1} \subset V_i \subset X
\]
and
\[
\dim V_i \leq \dim V_{i-1} + 1.
\]
Furthermore $\dim V_i = \dim V_{i-1}$ if and only 
if all the irreducible components
of $\overline{
\tilde{p}(\fibre{\tilde{q}}{\tilde{q}(\fibre{\tilde{p}}{V_{i-1}})}}$ 
have dimension $\dim V_{i-1}$.
Since these components all contain $V_{i-1}$ we see 
that $\dim V_i = \dim V_{i-1}$
if and only if 
\[
\overline{\tilde{p}(\fibre{\tilde{q}}{\tilde{q}(\fibre{\tilde{p}}{V_{i-1}})}}
=V_{i-1},
\]
{\em i.e.}\  if and only if $V_{i-1}$ is $\Hilb_x$-stable.
Since $\dim V_i \leq \dim X$ for all 
$i \in \{ 0, \ldots, n+1 \}$, we see that necessarily 
$\dim V_n = \dim V_{n+1}$,
so $V_n$ is $\Hilb_x$-stable. We set
\[
F_{x,y} := V_n.
\]
Let us now show that $F_{x,y}$ is minimal:
let $M$ be a $\Hilb_x$-stable subvariety through $y$. 
We will show inductively that $V_i \subset M$, the start
of the induction being clear since $y \in M$ by hypothesis. 
For the induction step, note that
if $V_{i-1} \subset M$, then 
\[
\overline{\tilde p(\fibre{\tilde q}{\tilde q(\fibre{\tilde p}{V_{i-1}})}} 
\subset \overline{\tilde p(\fibre{\tilde q}{\tilde q(
\fibre{\tilde p}{M})}}. 
\]
Since $M$ is $\Hilb_x$-stable the right hand side equals $M$, so 
\[
V_i \subset 
\overline{\tilde p(\fibre{\tilde q}{\tilde q(\fibre{\tilde p}{V_{i-1}})}} 
\]
implies the claim.

{\it Step 2: construction of the foliation.}
The preceding step gives the existence of a minimal 
$\Hilb_x$-stable subvariety 
$F_{x,y}$ through a general point $y \in X$.
Since the Chow scheme of $X$ has only countably 
many components, there exists a closed subvariety
$Z_x \subset \Chow X$ such that for $y \in X$ general, 
the variety $F_{x,y}$ corresponds to a point $[F_{x,y}] \in Z_x$.
Up to replacing $Z_x$ by a subvariety, we can suppose that the points 
$[F_{x,y}]$ are dense in $Z_x$.

Let $F_x \subset Z_x \times X$ be the universal 
family over $Z_x$. Note that since $Z_x$ is irreducible and 
the general fibre is irreducible, the
total space $F_x$ is irreducible. We claim that the natural morphism
\[
\holom{p'}{F_x}{X}
\]
is birational. We argue by contradiction and 
suppose that this is not the case: then for $y \in X$ general, the
fibre $\fibre{p'}{y}$ is not a singleton, so  $q(\fibre{p'}{y})$ 
has at least two distinct points $z$ and $z'$.
Thus by construction the subvarieties parametrised by $z$ and $z'$ 
are distinct minimal  $\Hilb_x$-stable
varieties through $y$, a contradiction to the uniqueness of these varieties.

Let $T_{F_x/Z_x}:= \Homsheaf(\Omega_{F_x/Z_x}, \sO_{F_x})$ 
be the relative tangent sheaf of the family. 
Since $p'$ is birational, the tangent map gives an integrable subsheaf  
\[
Tp' (T_{F_x/Z_x}) \subset T_X.
\]
We define the foliation $\sF_x$ to be the saturation 
of $Tp'(T_{F_x/Z_x})$ in $T_X$, {\em i.e.}\  the kernel of the surjective map
$T_X \rightarrow (T_X/Tp'(T_{F_x/Z_x}))/\mbox{Torsion}$.
By construction, the unique $\sF_x$-leaf  
through a general point $y \in X$ is the variety $F_{x,y}$, 
so it is the minimal $\Hilb_x$-stable subvariety. This 
also guarantees that the foliation is unique.
\end{proof}

\begin{remark}\label{remarkcano} 
In the preceding proof, we choose 
an irreducible component of 
$\overline{\tilde p(\fibre{\tilde q}{\tilde q(\fibre{\tilde p}{V_{i-1}})}}$ 
at each step. The resulting foliation
does not depend on these choices.
\end{remark}

\subsection{Examples and first properties}

It is of course possible that the foliation $\sF_x$ is trivial, that is $\sF_x=T_X$.
In this ``general type'' case our techniques don't say anything, but the following examples show that 
the geometrically most interesting cases are not of this type.

\begin{example}\label{examplerankone}
\begin{enumerate}
\item[a)] Let $X=\PP^n$ and consider the family of lines. 
Then for every $x \in \PP^n$, the foliation $\sF_x$ is the rank one 
foliation whose leaves are
the lines through $x$.

\item[b)] Let $X$ be projective manifold 
and let $l \subset X$ be a quasi-line such 
that $e(X,l)=1$. Then the morphism $p$ from Diagram \eqref{basicdiagram} 
is birational
and a general quasi-line meets the exceptional locus of $p$
in $x$ only. Thus the foliation $\sF_x$ has rank one, 
its general leaves are the general quasi-lines through $x$.
Conversely, if the foliation $\sF_x$ has rank one, the
leaf through a general point $y$ is a single curve, therefore 
$e(X,l)=1$.
\end{enumerate}
\end{example}

The following example was the starting point of our theory. It illustrates the philosophy that if we have a quasi-line
$l \subset X$ such that $e(X,l)>1$, the foliation $\sF_x$ should come from a fibration.

\begin{example} \label{examplepardini}
Let $X$ be a double cover of $\PP^1 \times \PP^2$
whose branch locus is a general divisor of bidegree $(2,2)$. 
The threefold $X$ is Fano
with Picard number $2$, 
denote by $\varphi~: X \to \PP^2$
the projection on the second factor. 
The map $\varphi$ is a conic bundle, whose
discriminant locus is a quartic curve 
in $\PP^2$.
Let $d$ be a general line in $\PP^2$ and set $S_d := \varphi ^{-1}(d)$.
The surface $S_d$ is a del Pezzo surface  
and the induced map $\varphi~: S_d \to d \simeq \PP^1$
has exactly $4$ singular fibres. 
Therefore we have a morphism \holom{\mu}{S_d}{\PP^2}
representing $S_d$ as the blow-up of
five points $p_1,\ldots,p_5 \in \PP^2$ in general position.
Let $C$ be a general line in $\PP^2$ and set $l := \mu^{-1}(C)$. 
Clearly $l$ is a quasi-line of $S_d$ and a section
of $\varphi~: S_d \to d$, therefore $l$ is a quasi-line of $X$ \cite[Lemma 4.1]{a3}.  
For every $x \in X$ general, the foliation $\sF_x$ is 
the rank two foliation whose general leaves are the preimages 
of lines through
$\varphi(x)$. Note that by consequence the singular locus of the 
foliation contains the conic $\fibre{\varphi}{\varphi(x)}$, 
but the general leaf is the smooth surface $S_d=\fibre{\varphi}{d}$.
\end{example}

As we already mentioned, the foliation $\sF_x$ is singular.
However we have the following properties.

\begin{proposition}\label{propositiontechnicalpropertiesx}
Let $x \in X$ be a general point 
and let $\sF_x$ be the foliation constructed in 
Theorem \ref{theoremexistencefoliation}.
Suppose that $\rk \sF_x < \dim X$.
Then the following properties hold.
\begin{enumerate}
\item[a)] A general quasi-line through $x$ meets 
the singular locus of $\sF_x$ exactly in $x$.
\item[b)] If $y \in X$ is a general point, all the quasi-lines 
parametrised by $\Hilb_x$ passing through $y$ are contained 
in $F_{x,y}$, 
the unique $\sF_x$-leaf through $y$. 
Furthermore they meet the singular locus of $\sF_x$ exactly in $x$.
\end{enumerate}
\end{proposition}

\begin{proof}
By construction the general $\sF_x$-leaf is 
of the form $F_{x,y}$, so it contains $x$. 
Since the foliation is not trivial, this implies 
that $x$ is in the singular locus $\sF_x^{\sing}$ of the foliation.
Since the foliation $\sF_x$ is saturated, its singular 
locus has codimension at least two. 
A nowadays well-known argument in deformation theory 
of very free rational curves \cite[II,Prop.3.7]{Kol96} shows that 
a general quasi-line $l$ through $x$ does not meet 
$\sF_x^{\sing} \setminus \{x\}$. 

For the second statement, the first part follows from
Remark \ref{remarkcano}.
Assume then by contradiction that for every $y \in X$ general 
there exists at least one
quasi-line through $x$ and $y$ that meets $\sF_x^{\sing} \setminus \{x\}$. 
We then get, with the notations of the Diagram \eqref{basicdiagram},
that for $y$ general in $X$, 
$$ p^{-1}(\sF_x^{\sing} \setminus \{x\}) \cap
q^{-1}(q(p^{-1}(y))) \neq \emptyset.$$
This means that a general quasi-line through $x$ intersects
$\sF_x^{\sing} \setminus \{x\}$, which is not possible since 
$\sF_x^{\sing} \setminus \{x\}$ is of codimension at least two
in $X$ (see again \cite[II,Prop.3.7]{Kol96}).
\end{proof}

We now give two technical consequences of our construction
which will be very useful in Section \ref{sectionmainresult}.

\begin{proposition} \label{propositionfactorisation}
With the notations above, the basic Diagram \eqref{basicdiagram} factors 
for $x$ general in $X$ through the universal family $F_x \rightarrow Z_x$ 
of the leaves of the foliation $\sF_x$, {\em i.e.}\  
up to replacing $\tilde{\Hilb}_x$ by a dense open subset
we have a commutative diagram 
\[
 \xymatrix{ 
\tilde{\Gamma}_x \ar[d]_{q|_{\tilde{\Gamma}_x}} \ar[r]^{p''} 
\ar @/^2pc/ [rr]^{p|_{\tilde{\Gamma}_x}} & F_x  \ar[d]^{q'} 
\ar[r]^{p'} & X \\
\tilde{\Hilb}_x  \ar[r]^{\psi} & Z_x & } 
\]
such that $\psi$ is dominant.
\end{proposition}

\begin{proof}
The natural morphism $p': F_x \rightarrow X$ is birational, 
so we have a rational map
\[
\merom{p'':=(p')^{-1} \circ p}{\Gamma_x}{F_x}.
\]
Since $\tilde{\Gamma}_x$ is smooth, we can replace $\Gamma_x$ 
by its normalisation without changing $\tilde{\Gamma}_x$.
Then \cite[Ch. 1.39]{Deb01} implies that the indeterminacy locus 
of $p''$ has codimension at least two, in particular it does not surject onto 
$\Hilb_x$. Therefore up to replacing $\tilde{\Hilb}_x$ by 
a dense open subset, 
we can suppose that $p''$ is defined on $\tilde{\Gamma}_x$. 
It is clear that $p|_{\tilde{\Gamma}_x}=p' \circ p''$,
and by construction of the foliation $\sF_x$, 
it maps a curve parametrised by $\tilde{\Hilb}_x$ 
into the general leaf containing it. So 
the rigidity lemma \cite[Lemma 1.15]{Deb01} applied to 
$\tilde{\Gamma}_x \rightarrow \tilde{\Hilb}_x$ 
and $\tilde{\Gamma}_x \rightarrow F_x \rightarrow Z_x$
shows that there exists a morphism 
$\psi: \tilde{\Hilb}_x \rightarrow Z_x$ that makes 
the diagram commutative. Since $p''$ is dominant, the same holds for $\psi$.
\end{proof}

Let us remark that the construction of 
the foliations $\sF_x$ can be done ``in family'' 
when $x$ moves in $X$. 
Indeed, for a general 
point $x \in X$, we have
defined a foliation $\sF_x$ and its universal family 
$q: F_x \rightarrow Z_x$. 
Arguing as in the construction of the foliation $\sF_x$, we see that
there exists a subvariety $Z \subset X \times \Chow X$
such that the projection on the first factor $\holom{p_X}{Z}{X}$
is a dominant proper morphism that satisfies $\fibre{p_X}{x}=Z_x$.
Thus we get a universal family $q:F \rightarrow Z$, where 
$$
F \subset Z \times X \subset X \times \Chow X \times X
$$
such that $\fibre{q'}{Z_x}=F_x$. The projection 
on the first and third factor induces a morphism
$\holom{p'}{F}{X \times X}$ such that the restriction 
$\holom{p'|_{F_x}}{F_x}{x \times X}$
identifies to the morphism $\holom{p'}{F_x}{X}$ 
from Diagram \eqref{basicdiagramfoliation}\footnote{In order 
to simplify the notation, we denote
the two morphisms by the same letter.}. 
The morphism $\holom{p'}{F}{X \times X}$ is surjective and birational.
Last but not least, let $\holom{p_{\mathcal{C}}}{Z}{\chow{X}}$ 
be the morphism given by the projection on the second factor. 
We resume the construction in a commutative diagram
\begin{equation} \label{basicdiagramglobal}
 \xymatrix{ 
& F \ar[d]_{q'}\ar[r]^{p'} & X \times X \ar[ldd]^{p_1} \\
\chow{X} & Z  \ar[d]_{p_X}\ar[l]^{p_{\mathcal{C}}} & \\
& X
}
\end{equation}

\section{A fundamental dichotomy}

If the foliation $\sF_x$ is not trivial, it is singular at $x$
and depends on the choice of the general base point $x$. 
We will show now that it is very interesting to compare 
the foliations arising from different general choices of base points.
The aim of this section is to prove the following dichotomy.

\begin{proposition}\label{propositiontwofoliations}
Let $X$ be a projective manifold containing
a quasi-line $l$.
Let $x$ and $y$ be two general points in $X$ and 
denote by $\sF_x$ and $\sF_y$ the corresponding foliations.
Assume $\sF_x$ and $\sF_y$ are not trivial, that is $\rk \sF_x=\rk \sF_y<\dim X$. 
Denote by $F_{x,y}$ the unique  $\sF_x$-leaf through $y$
and by $F_{y,x}$ the unique $\sF_y$-leaf through $x$.

Then the following properties hold.
\begin{enumerate}
\item[a)] The foliation $\sF_x$ is smooth 
at $y$.
\item[b)] If $y$ is very general, 
there exists a desingularisation   
$\mu : F_{x,y}' \rightarrow F_{x,y}$ such that $y$ is in 
the locus of free rational curves (cf. \cite[Prop.4.14]{Deb01}) of $F'_{x,y}$.
\item[c)] The foliation $\sF_y$ is smooth at $x$.
\end{enumerate}

Moreover {\em exactly one} of the following situations occurs.
\begin{enumerate}
\item[{\em (I)}] We have $F_{x,y}=F_{y,x}$ and 
the variety $F_{x,y}$ is smooth at $x$.
\item[{\em (II)}] The intersection $F_{x,y} \cap F_{y,x}$ 
is a strict subset of $F_{x,y}$ and $F_{y,x}$.
\end{enumerate}
\end{proposition}

\begin{proof}
It is obvious how to ensure the properties a) and b): 
if $F_{x,z}$ is a general $\sF_x$-leaf,
take a desingularisation $\mu : F_{x,z}' \rightarrow F_{x,z}$. 
Then a very general point $y \in F_{x,z}'$ is in the free locus and does 
not meet the exceptional locus of $\mu$, so it can 
be considered as a very general point of $y \in F_{x,z}$.
Since for such a very general point $F_{x,y}=F_{x,z}$, we get the conclusion.

Let us now show that for $y \in X$ general, property c) also holds:
let $R \subset X \times X$ be the smallest closed subset 
such that for $y \in X$ general, the fibre 
$R_y:=\fibre{p_1}{y} \subset y \times X$ induced 
by the projection on the first factor is the 
singular locus of the foliation $\sF_y$. 
Since the foliations $\sF_y$ are saturated, 
we have $\dim R \leq 2 \dim X - 2$. 
Let \holom{p_2}{R}{X} be the restriction of the projection 
on the second factor. For $x \in X$ general, the fibre
$\fibre{p_2}{x}$ has dimension at most $\dim X-2$. Therefore 
$p_1(\fibre{p_2}{x})$ has dimension 
at most $\dim X-2$, so it is a strict subset of $X$. Yet by construction,
$$
\{ y \in X \ \mbox{general}  \ | \ \sF_y \ 
\mbox{is singular at} \ x \} \subset p_1(\fibre{p_2}{x}).
$$

Let us now prove the basic dichotomy: 
suppose that the intersection $F_{x,y} \cap F_{y,x}$ 
is not a strict subset of  $F_{x,y}$ or $F_{y,x}$. 
Since the foliations $\sF_x$ and $\sF_y$ have the same rank, 
this implies $F_{x,y}=F_{y,x}$.
Furthermore we have just seen that the point $x$ is in the regular
locus of the foliation $\sF_y$. 
Thus $F_{y,x}$ is smooth at $x$, hence $F_{x,y}=F_{y,x}$ is smooth at $x$.
\end{proof}

\begin{examples}
\begin{enumerate}
\item[a)] We don't have any examples for the Case (II) 
of Proposition \ref{propositiontwofoliations}. 
Nevertheless we see no reason why such examples
shouldn't exist, since the foliations depend heavily 
on the family of curves and thus on the choice of the base point.
\item[b)] Let $X=\PP^n$ and $\Hilb$ be the family of lines. 
Then for every $x,y \in \PP^n$ with $x \neq y$, the 
varieties $F_{x,y}=F_{y,x}$ are the unique line through $x$ and $y$.
More generally, case (I) occurs when $e(X,l)=1$.
\item[c)] In Example \ref{examplepardini}, 
for 
$x \in X$ general, the general $\sF_x$-leaf is a smooth surface 
$S_d=\fibre{\varphi}{d}$ where $d$ is a line in $\PP^2$ through $\varphi(x)$.
Thus for $x$ and $y$ in $X$ general, 
the varieties $F_{x,y}$ and $F_{y,x}$ are 
the preimage \fibre{\varphi}{d_{\varphi(x),\varphi(y)}} of the 
unique line $d_{\varphi(x),\varphi(y)}$ through $\varphi(x)$ and 
$\varphi(y)$.
\end{enumerate}
\end{examples}

\begin{remark} \label{remarksecondcaseglobal}
Using the Diagram \eqref{basicdiagramglobal}, 
we can give a technically more useful expression of the dichotomy
in Proposition \ref{propositiontwofoliations}: 
fix $[\sL]$ a general point in $p_{\mathcal{C}}(Z)$
(which also means that $\sL$ is a general 
$\sF_x$-leaf, $x$ being general).
By definition of the foliations and the parameter space $Z$, we have
\[
\fibre{p_{\mathcal{C}}}{[\sL]} = 
\overline{\{ (z,[\sL]) 
\in Z  \ | \ \sL \ \mbox{is a} \ \sF_z-\mbox{leaf} \}}.
\]
Being in case (I) of Proposition \ref{propositiontwofoliations} 
is equivalent to 
have $F_{x,z}=F_{z,x}$ for $x$ and $z$ general in $X$.
Since for $z \in \sL$ general, we have $F_{x,z}=\sL$, 
we see that 
being in case (I) of Proposition \ref{propositiontwofoliations} is equivalent 
to asking that 
\[
p_X (\fibre{p_{\mathcal{C}}}{[\sL]})= \sL.
\]
Since $p_X$ is an isomorphism on the fibres of 
$p_{\mathcal{C}}$, this is equivalent to
\[
\fibre{p_{\mathcal{C}}}{[\sL]} = \sL \times [\sL].
\]
Looking at the universal family, it is clear that
\[
\fibre{q'}{\fibre{p_{\mathcal{C}}}{[\sL]}} =
\fibre{p_{\mathcal{C}}}{[\sL]}\times \sL,
\]
so we get a natural identification
\[
\fibre{q'}{\fibre{p_{\mathcal{C}}}{[\sL]}} = 
\sL \times [\sL] \times \sL.
\]
This means that $\sL$ is the natural parameter space for the $x \in X$ such that 
the variety $\sL \subset X$ is a leaf of the foliation $\sF_x$.  
\end{remark}

\section{The main results}
\label{sectionmainresult}

The aim of this section is to prove Theorem \ref{maincaseI} and Proposition \ref{propositioncorankone}.
Until the end of Subsection \ref{subsectionproofs}, we will always suppose that the Assumption \ref{assumption}
holds, {\em i.e.}\  we are in the following situation.

{\bf Assumption.}
Let $X$ be a projective manifold containing
a quasi-line $l$.
Let $x$ and $y$ be two general points in $X$ and 
denote by $\sF_x$ and $\sF_y$ the corresponding foliations.
Denote by $F_{x,y}$ the unique  $\sF_x$-leaf through $y$
and by $F_{y,x}$ the unique $\sF_y$-leaf through $x$.
Then  we have 
$$
F_{x,y}=F_{y,x}.
$$ 
We also assume that the foliation $\sF_x$ is not trivial ($\rk \sF_x<\dim X$),
and denote by $\sL$ a general $\sF_x$-leaf.
Since our assumption means precisely that we are in the first case of 
Proposition \ref{propositiontwofoliations}, the leaf $\sL$ is smooth in $x$.

\subsection{Structure of the general leaves.}

Recall the commutative diagram 
\begin{equation}
\label{basicdiagramfoliation}
\xymatrix{ 
\tilde{\Gamma}_x \ar[d]_{q|_{\tilde{\Gamma}_x}} \ar[r]^{p''} 
\ar @/^2pc/ [rr]^{p|_{\tilde{\Gamma}_x}} & F_x  \ar[d]^{q'} 
\ar[r]^{p'} & X \\
\tilde{\Hilb}_x  \ar[r]^{\psi} & Z_x & } 
\end{equation}
introduced in Proposition \ref{propositionfactorisation}.

\begin{lemma}
\label{lemmasecondcasequasilines}
Under the Assumption \ref{assumption} the following holds. 
\begin{enumerate}
\item[a)] A general quasi-line through $x$ contained in $\sL$ 
is contained in the smooth locus of $\sL$ 
and is a quasi-line in that variety.
\item[b)] In particular if $y$ is a general point in $\sL$, every quasi-line through $x$ and $y$ is contained 
in the smooth locus of $\sL$ and is a quasi-line in that variety.
\end{enumerate}
\end{lemma}

\begin{remark*}
It is probably necessary to explain this statement: 
the leaf $\sL$ corresponds to a general point $[\sL] \in Z_x$. Let 
$\psi: \tilde{\Hilb}_x \rightarrow \tilde{Z}_x$ be the morphism in the 
diagram \eqref{basicdiagramfoliation} above, then the fibre
$\fibre{\psi}{[\sL]}$ parametrises quasi-lines through $x$ 
contained in $\sL$. The term general  
(see Definition \ref{definitiongeneral})
refers to a quasi-line parametrised by a general 
point of $\fibre{\psi}{[\sL]}$. Since $\sL$ 
is a general leaf, this is equivalent
to saying that the quasi-line corresponds to 
a general point of $\tilde{\Hilb}_x$.

This also explains why the second statement is a special case of the first: 
the choice of the point $x$ and the $\sF_x$-leaf $\sL$ fix $\fibre{\psi}{[\sL]}$,
but the choice of a general $y \in \sL$ is independent, so the  
quasi-lines through $x$ and $y$ are general.
\end{remark*}

\begin{proof}
By Proposition \ref{propositiontechnicalpropertiesx} (a), 
a general quasi-line $l$ through $x$ meets $\sF_x^{\sing}$ exactly in $x$. 
Since 
\[
\sL^{\sing} \subset \sF_x^{\sing}
\]
and $x \not\in \sL^{\sing}$ by Proposition \ref{propositiontwofoliations}, the curve $l$
is contained in the smooth locus of $\sL$.

Since $[\sL] \in Z_x$ is a general point, 
all the irreducible 
components of the fibre $\fibre{\psi}{[\sL]}$
have the expected dimension $\rk \sF_x -1$. 
Let $V$ be such an irreducible component, then \fibre{q}{V} 
has dimension $\rk \sF_x$.
Since $[\sL]$ is general, the variety $\fibre{q}{V}$ is not 
contained in the exceptional locus of the generically finite morphism $p$, so
$p(\fibre{q}{V})$ has dimension $\rk \sF_x$. Since we have a 
factorisation $p=p' \circ p''$, we see that $p''(\fibre{q}{V})$ 
has dimension $\rk \sF_x=\dim \sL$.
So the family of rational curves parametrised by $V$ dominates 
$\sL$ and all its members pass through $x$. Thus
by \cite[4.10]{Deb01} a curve parametrised by a general point 
of $V$ is a very free curve in $\sL$.
Since $p'|_\sL$ is an isomorphism, this implies that 
$T_\sL|_l$ is ample and the injection
\[
T_\sL|_l \hookrightarrow T_{X}|_l
\]
implies that it has the splitting type of a quasi-line. 
\end{proof}

Our goal is now to show the following proposition which corresponds to part a) of Theorem \ref{maincaseI}.

\begin{proposition} \label{propositionleaffibration}
Under the Assumption \ref{assumption}, suppose moreover that
$$
\dim_x \sF_x^{\sing} < \rk \sF_x-1+ \frac{1}{2} (n-\rk \sF_x).
$$
Then there exists a rational fibration \merom{\varphi_\sL}{(\sL, l)}{\PP^1} such that $l$ is a section.
\end{proposition}

The statement will be a consequence of Theorem \ref{theoremvoica}. In view of the hypothesis in 
Theorem \ref{theoremvoica} we have to address the following tasks:

\begin{enumerate}
\item construct an effective divisor $E_{\sL,x} \subset \sL$ such that $E_{\sL,x} \cdot l=1$,
\item show that $E_{\sL,x}$ is in a linear system of dimension one.
\end{enumerate}

We will see that $E_{\sL,x}$ will be the restriction of an exceptional divisor $E_x$ of the natural map 
$F_x \rightarrow X$. At first glance it seems impossible that  $E_{\sL,x}$ moves in a linear system, but this
intuition is wrong: the restriction of $F_x \rightarrow X$ to $\sL$ is an embedding, so $E_{\sL,x}$ is not exceptional.

If the foliation $\sF_x$ has rank one, the construction of $E_x$ is trivial: the universal families $F_x$ and $\Gamma_x$ are
the same, the set $\fibre{p}{x} \subset \Gamma_x$ gives a $q$-section, so its restriction to a general quasi-line is a reduced point. If $\sF_x$ has higher rank, this obvious construction no longer works:
the generically finite morphism $p''$ maps  $\fibre{p}{x}$ onto $\fibre{p'}{x}$ which is  not a divisor in $F_x$.

{\em Step 1. Construction of the divisor $E_{\sL,x}$.}
We use the notation of the Diagram \eqref{basicdiagramfoliation}.
Any $\sF_x$-leaf passes through $x$, so $\fibre{p'}{x}$ is a section of $q'$.
Furthermore $\sL$ is smooth at $x$, so the general $q'$-fibre is smooth at the intersection with $\fibre{p'}{x}$.
This shows that $F_x$ is smooth in the general point of $\fibre{p'}{x}$. 
Since by Lemma \ref{lemmasecondcasequasilines} above a general curve parametrised by $\tilde{\Hilb}_x$ is
contained in the smooth locus of a general $q'$-fibre, it follows that it is contained in the smooth locus of $F_x$.
This shows that there exists a desingularisation
$F'_x \rightarrow F_x$ that 
\begin{enumerate}
\item[a)] is an isomorphism in the general point of \fibre{p'}{x}, 
\item[b)] a general curve parametrised by $\tilde{\Hilb}_x$ does not meet the exceptional locus.
\end{enumerate} 
These properties assure that all the following computations take place in the complement of the exceptional locus
of the desingularisation, moreover we don't have to distinguish Weil and Cartier divisors. 
In order to simplify the notation, we thus suppose without loss of generality that $F_x$ is smooth.

Since the base of the birational morphism $\holom{p'}{F_x}{X}$ is smooth, the $p'$-exceptional locus
has pure codimension one in $F_x$ and the determinant of the Jacobian $d p'$ defines an equality of cycles 
\[
K_{F_x/X} = \sum_{i \geq 1} a_i E_i,
\]
where the $E_i$ are the $p'$-exceptional divisors and the $a_i$ are positive integers.
Let $l \subset X$ be a quasi-line parametrised by a general point of $\tilde{\Hilb_x}$ such that $l \subset \sL$. 
Identifying the leaf $\sL \subset X$ and the corresponding subvariety of $F_x$, 
we authorise ourselves to see $l$ also as a rational curve in $F_x$.
Then $l$ is a free curve in $F_x$ (although it is not a quasi-line in $F_x$), so 
\[
E_i \cdot l \geq 0 \qquad \forall \ i \geq 1.
\] 
Note furthermore that 
\[
p'(E_i) \subset \sF_x^{\sing},
\]
since through any point of $p'(E_i)$ pass an infinity of leaves.
By Proposition \ref{propositiontechnicalpropertiesx} a general quasi-line through $x$ meets $\sF_x^{\sing}$ exactly in $x$, so we see that
\[
E_i \cdot l > 0 \ \Leftrightarrow \   \fibre{p'}{x} \subset E_i.
\] 
Thus we get
\[
K_{F_x/X} \cdot l = \sum_{\fibre{p'}{x} \subset E_i} a_i E_i \cdot l
\]
and we denote by $E_x$ the support of $\sum_{\fibre{p'}{x} \subset E_i} a_i E_i$. 
Note that $E_x$ is not empty, since 
$x$ is contained in every $\sF_x$-leaf.
Furthermore every irreducible component $E_i \subset E_x$ maps surjectively onto $Z_x$ since 
$\fibre{p'}{x}$ is a $q'$-section. 
In particular the general fibre of $E_i \rightarrow Z_x$ has dimension $\rk \sF_x-1$, so 
$$
E_{\sL,x} := E_x \cap \sL
$$ 
is a non-empty effective divisor in $\sL$ passing through the point
$x \in \sL$.
Since the restriction of $p'$ to $\sL$ is an isomorphism, we see that 
\[ 
\dim p'(E_i) \geq \rk \sF_x-1,
\]
so in particular 
\[
\dim_x \sF_x^{\sing} \geq \rk \sF_x-1.
\]
The following lemma shows that the divisor $E_{\sL,x}$ has some interesting intersection properties
provided the dimension of the singular locus of the foliation is not too large.

\begin{lemma} \label{lemmaalmostline} Under the Assumption \ref{assumption},
let $l$ be a quasi-line parametrised by a general point of $\tilde{\Hilb_x}$ such that $l \subset \sL$. 
If
\[
\dim_x \sF_x^{\sing} < \rk \sF_x-1+ \frac{d-1}{d} (n-\rk \sF_x)
\]
for some $d \in \N$, then
\[
E_{\sL, x} \cdot l < d.
\]
In particular if 
$$
\dim_x \sF_x^{\sing} < \rk \sF_x-1+ \frac{1}{2} (n-\rk \sF_x),
$$
then 
\[
E_{\sL, x} \cdot l=1.
\]
Hence  the effective divisor $E_{\sL, x}$ is smooth in $x \in \sL$ and $E_x$ is irreducible. 
\end{lemma}

\begin{proof}
By Lemma \ref{lemmasecondcasequasilines} the rational curve 
$l$ is a quasi-line in the leaf $\sL$, so by adjunction
\[
K_{F_x} \cdot l = K_{\sL} \cdot l = - (\rk \sF_x+1).
\]
Since $l$ is a quasi-line in $X$ this implies
\[
K_{F_x/X} \cdot l = n-\rk \sF_x.
\]
Using the notation introduced before the lemma, a local computation shows that
if $E_i \subset E_x$ is an irreducible component then  
\[
a_i \geq n-1-\dim p'(E_i).
\]
Since $x \in p'(E_i) \subset \sF_x^{\sing}$ we get
\[
a_i \geq n-1-\dim_x \sF_x^{\sing}.
\]
Thus one has
\[
n-\rk \sF_x =  K_{F_x/X} \cdot l = \sum_{\fibre{p'}{x} \subset E_i} a_i E_i \cdot l \geq (n-1-\dim_x \sF_x^{\sing}) E_{\sL, x} \cdot l.
\]
So if $E_{\sL, x} \cdot l \geq d$, we have
\[
n-\rk \sF_x \geq d (n-1-\dim_x \sF_x^{\sing}).
\]
This is equivalent to the first statement.

The second statement corresponds to the case $d=2$, the smoothness
of $E_{\sL, x}$ in $x$ is immediate from $E_{\sL, x} \cdot l=1$. Since any irreducible component of $E_x$ contains
$\fibre{p'}{x}$, this also shows that $E_x$ is irreducible.
\end{proof}

{\em Step 2. Dimension of the linear system.}
Recall the commutative Diagram \eqref{basicdiagramglobal}
\[
 \xymatrix{ 
& F \ar[d]_{q'}\ar[r]^{p'} & X \times X \ar[ldd]^{p_1} \\
\chow{X} & Z  \ar[d]_{p_X}\ar[l]^{p_{\mathcal{C}}} & \\
& X
}
\]
The map $p'$ is birational and since 
$X \times X$ is smooth, the $p'$-exceptional locus has codimension 1
in $F$. Denote by
$E \subset F$ the union of $p'$-exceptional divisors that 
contain $\fibre{p'}{\Delta}$. By definition we have for $x \in X$ general
\[
E_x = \fibre{(p_X \circ q')}{x} \cap E,
\]
where $E_x \subset F_x$ is the divisor introduced in the first step.
By Remark \ref{remarksecondcaseglobal} one has
\[
\fibre{q'}{\fibre{p_{\mathcal{C}}}{[\sL]}} = 
\sL \times [\sL] \times \sL,
\]
so $\sL$ is the natural parameter space for the $x \in X$ such that the subvariety
$\sL \subset X$ is a leaf of the foliation $\sF_x$.  
Thus for $x$ general in $\sL$ we have
$$ 
E_{\sL,x}
= 
\{ y \in \sL \mid (x,[\sL],y) \in E\}
$$
where $E_{\sL, x} \subset \sL$ is the divisor introduced in the first step.
Varying $x$ in the parameter space $\sL$, this defines a family of algebraically equivalent divisors on the leaf 
$\sL$,
and since $x \in E_{\sL,x}$
the divisor $E_{\sL,x}$ moves in $\sL$.
Since $\sL$ is rationally connected, we have $h^1(\sL,
\mathcal O _{\sL})=0$, 
the algebraically equivalent divisors $E_{\sL,x}$ 
are therefore linearly equivalent.
Simply denoting by $E_{\sL}$ any linear representative of
the $E_{\sL,x}$'s, we finally deduce that 
\[
h^0(\sL, E_{\sL}) \geq 2.
\]

\begin{proof}[End of the proof of Proposition \ref{propositionleaffibration}] 
According to Theorem \ref{theoremvoica}, the
only thing to prove is the equality $h^0(\sL, E_{\sL}) = 2$.
Set $h^0(\sL, E_{\sL}) =: s+1$,
and let $\varphi _{\sL} : \tilde{\sL} \to \mathbb P ^s$
be the map given by Theorem \ref{theoremvoica}.
Recall now the construction of $\sL$ by a sequence of subvarieties 
 \[
V_0=y \subset V_1 \subset \ldots \subset V_{n}= 
\sL \subset X
\]
where $V_i$ is obtained form $V_{i-1}$ by adding quasi-lines 
through $x$ and a general point of $V_{i-1}$. We will prove inductively that
$\varphi _{\sL}(V_i)$ is a fixed line in $\mathbb P ^s$
for every $i \in \{ 1, \ldots, n \}$.
Since $V_1$ is a quasi-line $l'$ through $x$ and a given point
of $\sL$, the start 
of the induction is trivial.
Suppose now that $\varphi _{\sL}(V_{i-1})=\varphi _{\sL}
(l')$ for some 
$i \geq 2$. The quasi-lines $l_z$ through $x$ and a 
general point $z$ of $V_{i-1}$ dominate
$V_i$, so it is sufficient to show that $\varphi _{\sL}(l_z)=
\varphi _{\sL}(l')$.
Yet $\varphi _{\sL}(l_z)$ and $\varphi _{\sL}(l')$ are 
lines in $\PP^s$ meeting in the points $\varphi _{\sL}(x)$ and 
$\varphi _{\sL}(z)$, thus they are identical. 
\end{proof}

In the preceding proof the hypothesis on the singular locus of the foliation was only
used to apply Lemma \ref{lemmaalmostline}. Thus we have shown:

\begin{proposition} 
Under the Assumption \ref{assumption}, suppose moreover 
that $E_{\sL}\cdot l =1$.
Then there exists a rational fibration \merom{\varphi_\sL}{(\sL, l)}{\PP^1} such that $l$ is a section.
\end{proposition}

\subsection{Structure of $X$.} \label{subsectionproofs}

The first part of Proposition \ref{propositioncorankone} is a corollary of Theorem \ref{maincaseI}: the only thing we have to show is that
the surface $Y$ admits a birational map $Y \dashrightarrow \PP^2$ such that $l'$ maps onto a line. This can be deduced from
the list in \cite[Prop. 1.21]{IV03}. Since the proof of Theorem \ref{maincaseI} is rather involved, we give a short independent argument
based on Theorem \ref{theoremvoica}.

\begin{proof}[Proof of Proposition \ref{propositioncorankone}]
Denote by $\sL$ a general leaf.
By Lemma \ref{lemmasecondcasequasilines} a general quasi-line $l$ contained in  $\sL$ is a quasi-line in $\sL$. Therefore the exact sequence
\[
0 \rightarrow T_{\sL}|_l \rightarrow T_X|_l \rightarrow N_{\sL/X}|_l \simeq \sO_X(\sL)|_l \rightarrow 0
\]
shows that $\sL \cdot l=1$. Since the divisor $\sL$ moves in the rationally connected manifold $X$, we can apply Theorem \ref{theoremvoica}
to get a  birational map $\sigma : \tilde X \to X$  and a 
fibration $\varphi : \tilde X \to \mathbb P^k$ with the stated properties, so the only thing to show
is $k=2$. Since $\varphi$ is induced by the moving part of the linear system $|\sigma^* \sL |$ the general leaf maps onto a hyperplane.
Furthermore it is not hard to see that the restriction of $|\sL|$ 
to a leaf is the linear system $|E_{\sL}|$ from Proposition \ref{propositionleaffibration}, 
so the restriction of $\varphi$ to $\sL$ is the fibration $\varphi_L$. Thus the hyperplane $\varphi(\sL)$ is a curve.

If the Picard number of $X$ is one, the effective divisor $\sL$ is ample.
By \cite[Thm.4.4, Cor.4.6]{BBI00} this implies that $X$ is a projective space and the quasi-lines are lines. 
The foliation defined by lines in the projective space has rank one, so $X \simeq \PP^2$.
\end{proof}

\begin{proof}[Proof of Theorem \ref{maincaseI}]
Statement a) of the theorem is settled by Proposition \ref{propositionleaffibration}. Suppose now that we are in
the situation of the Statement b).
Fix a general point $x \in X$ and consider the basic diagram
$$
\xymatrix{ 
F_x  \ar[d]_{q'}  \ar[r]^{p'} & X \\
Z_x & } 
$$
where $F_x$ is the universal family of the foliation $\sF_x$.
We have shown in Proposition \ref{propositionleaffibration} 
that for $[\sL]$ in $Z_x$ general, 
there exists a complete linear system $|E_\sL|$
of dimension one. Since we have fixed $x \in X$, we also have a distinguished element $[E_{\sL, x}] \in |E_\sL|$,
where
$$
E_{\sL, x} := E_{x} \cap \sL
$$
is the family of cyles introduced in Lemma \ref{lemmaalmostline}.

Since $h^1(\sL, \sO_\sL)=0$ 
the line bundle $E_\sL$ does not deform, so $|E_\sL|$
is an irreducible component of the Chow scheme $\chow{\sL}$. 
By countability of the number of irreducible components
of the relative Chow scheme $\chow{F_x/Z_x}$ 
there exists an integral scheme $U_x \rightarrow Z_x$ such that
the general fibre identifies to $|E_\sL| \simeq \PP^1$. 
The fibration $U_x \rightarrow Z_x$ has a distinguished (rational) section $S_x$ given by $x \mapsto [E_{\sL,x}] $.  

Consider now $U_x$ as a subset of $\chow{F_x}$. Then by \cite[I. Thm.6.8]{Kol96} the holomorphic map $p'$ induces a map 
\holom{p'_*}{\chow{F_x}}{\chow{X}} and we denote by $Y_x$ the image of $U_x$. It is immediate that $U_x \rightarrow Y_x$ is birational.
Let us now look at the family of cycles $E_{\sL, x}$: 
clearly their image by $p'$ is contained 
in the variety $p'(E_x)$ which has dimension at least $\rk \sF_x-1$,
since the restriction of $p'$ to a general leaf $E_{\sL, x}$ is an isomorphism.
Since
$p'(E_x) \subset \sF_x^{\sing}$ the hypothesis on the dimension of the singular locus implies that $p'(E_x)$ 
has dimension exactly $\rk \sF_x-1$. Thus the image of $E_{\sL, x}$ 
does not depend on $\sL$. This shows that the distinguished section $S_x$ is contracted by $p'_*$. 

Let $Y$ be the normalisation of $Y_x$ and denote by $\tilde{X}$ the normalisation of the universal family over $Y$.
Denote by \holom{\sigma}{\tilde{X}}{X} and \holom{\varphi}{\tilde{X}}{Y} the natural maps so that we get  a commutative diagram
 $$
\xymatrix{ 
F_x  \ar[d]_{q'}  \ar[r]^{p'} & X  &  \tilde{X} \ar[l]_\sigma \ar[d]^\varphi  \\
Z_x \ar[r]_{p'_*} &  Y_x & Y \ar[l] } 
$$
Since $q'(\fibre{p'}{x}) \subset S_x$ and the normalisation $Y \rightarrow Y_x$ is finite, a diagram chase shows that
$\fibre{\sigma}{x}$ is finite.
Yet $\sigma$ is birational and $X$ is normal,
so it follows by Zariski's main theorem that $\fibre{\sigma}{x}$ is a singleton and $\sigma$ is an isomorphism in a 
neighbourhood of $x$. In order to simplify the notation, we identify $\fibre{\sigma}{x}$ and $x$.
By construction the fibre $\fibre{\varphi}{\varphi(x)}$ is equal to $E_{\sL, x}$, where $\sL$ is a general leaf.
Since $\tilde{X}$ is smooth in $x$ and  $E_{\sL, x}$ is smooth in $x$ by Lemma \ref{lemmaalmostline}, the normal variety $Y$
is smooth in $\varphi(x)$.

Let now $l \subset \tilde{X}$ be a general quasi-line through $x$ contained in a general leaf $\sL$.
Since the birational map $\sigma$ is an isomorphism around $x$, the quasi-line $l$
does not meet the locus where $\sigma^{-1}$ is not a morphism.
Thus we can see $l$ as quasi-line in $\tilde{X}$ passing through $x$ and we will now show that
$l':=\varphi(l)$ is a quasi-line in $Y$. Note first that since $Y$ is smooth in $\varphi(x)$, the rational curve
$l'$ is contained in the smooth locus of $Y$.
The restriction of $\varphi$ to 
$\fibre{\varphi}{\varphi(l)}$ identifies to the graph of the meromorphic fibration
$\varphi_\sL: \sL \dashrightarrow \PP^1$ constructed in Proposition \ref{propositionleaffibration}. 
In particular $l$ is a quasi-line in $\fibre{\varphi}{\varphi(l)}$, so $l' \subset Y$ is a quasi-line
by \cite[Lemma 4.1]{a3}.

Let $\Hilb'_{\varphi(x)}$ be the normalisation of the subset of the Chow scheme $\chow{Y}$
parametrising deformations and degenerations of $l'$ passing through $y$. The push-forward induces a rational map
\merom{\varphi_*}{\Hilb_x}{\Hilb'_{\varphi(x)}} and a dimension count (see the proof of \cite[Lemma 4.1]{a3})
shows that this map is dominant. 
Apply now Theorem \ref{theoremexistencefoliation} to the family of quasi-lines $l' \subset Y$ through $\varphi(x)$
and denote by $\sG_{\varphi(x)}$ the resulting foliation on $Y$.
Since the general leaves of the foliations $\sF_{x}$  ({\em resp.} $\sG_{\varphi(x)}$) are built
from quasi-lines parametrised by $\Hilb_x$ ({\em resp.} $\Hilb'_{\varphi(x)}$),
we can use the dominant map \merom{\varphi_*}{\Hilb_x}{\Hilb'_{\varphi(x)}} to show that
a general $\sF_x$-leaf is the $\varphi$-preimage of a general $\sG_{\varphi(x)}$-leaf (the tedious details are left to the reader).
Since $Y$ has dimension $n-\rk \sF_x+1$, this implies that the general $\sG_{\varphi(x)}$-leaf has rank one.
Thus $e(Y,l')=1$ by Example \ref{examplerankone}.
\end{proof}

\subsection{Examples}

The following example generalises Example \ref{examplepardini}
to a situation where the quasi-line $l' \subset Y$ is not necessarily a line
in a projective space.

\begin{example} \label{exampleconicbundle}
Let $Y$ be a projective manifold of dimension $n-1$ 
and let $l' \subset Y$ be a quasi-line such that $e(Y,l')=1$.
Let now $\holom{\varphi}{X}{Y}$ be a conic bundle 
over $Y$ such that the discriminant locus $\Delta \subset Y$ 
satisfies $\Delta \cdot l' \geq 2$. For a general quasi-line $l'$, 
the preimage $X_{l'}:=\fibre{\varphi}{l'}$ is a smooth surface.  
Thus by \cite[Lemma 4.1]{a3}, there exists a quasi-line
$l \subset X$ that is a section of  $X_{l'} \rightarrow l'$.
Let $x$ be a general point of $X$.
Since $e(Y,l')=1$, the $\sF_x$-leaves 
are contained in the surfaces  
$\fibre{\varphi}{l'}$, hence 
$\rk \sF_x \leq 2$.
Moreover if we had $e(X,l)=1$ an argument analogous to the proof of \cite[Thm.4.14]{a3}
would show that $\holom{\varphi}{X}{Y}$ is smooth. Since this is not the case, we have $e(X,l)>1$,
so the foliation $\sF_x$ has rank two.
Note also that if $\sG_{\varphi(x)}$ denotes the rank one foliation associated to $l' \subset Y$, then
$$
\dim_x \sF_x^{\sing} = \dim_{\varphi(x)} \sG_{\varphi(x)}^{\sing} + 1.
$$
\end{example}

We will now construct an example of a fourfold $X$ containing a quasi-line $l$ such that
 for $x \in X$ general one has
$\rk \sF_x=1$ and $\dim_x \sF_x^{\sing}=1$, hence
\[
\dim_x \sF_x^{\sing} < \rk \sF_x-1+ \frac{1}{2} (n-\rk \sF_x),
\]
but $\dim_x \sF_x^{\sing} \neq \rk \sF_x-1$. Using the construction in Example \ref{exampleconicbundle}
one can transform this into an example where the foliation has rank two.

\begin{example}
Let $d_0$ be a line in $\PP^3$ and let $E$ be a rank two vector bundle given by the Serre extension
\[
0 \rightarrow \sO_{\PP^3} \rightarrow E \rightarrow \sI_{d_0}(-1) \rightarrow 0.
\]
Set $X:=\PP(E)$, denote by $\holom{\varphi}{X}{\PP^3}$ the natural projection, and by $\xi$
the unique effective divisor in the linear system $|\sO_{\PP(E)}(1)|$. By the canonical bundle formula
\[
-K_X = \varphi^* 5H + 2 \xi,
\]
and one checks easily that $-K_X$ is nef.

If $d \subset \PP^3$ is a general line, then $E|_d \simeq \sO_{\PP^1} \oplus \sO_{\PP^1}(-1)$,
so by \cite[Prop.5.1.,Rem.5.6.]{a3} the variety $X$ contains a quasi-line $l$ such that
$\varphi(l)$ is a line and $e(X,l)=1$. Thus if we fix a general point $x \in X$ (in particular $x \notin \xi$), the foliation 
$\sF_x$ has rank one and we claim that 
$$
\sF_x^{\sing} = \fibre{\varphi}{\varphi(x)} \cup (\fibre{\varphi}{H_x} \cap \xi),
$$ 
where $H_x \subset \PP^3$ is the unique hyperplane containing $\varphi(x)$ and $d_0$.
In particular $\sF_x^{\sing}$ has dimension one in $x$. 

{\em Proof of the claim.} 
We start by analysing the degenerations of $l$ passing through $x$:
we have $-K_X \cdot l=5, \varphi^* H \cdot l=1$, so any degeneration 
is of the form 
\[
l_0 + \sum_{i \geq 1} l_i
\]
such that $\varphi^* H \cdot l_0=1$ and $\varphi^* H \cdot l_i=0$ for all $i>0$. Thus for $i>1$ the curve $l_i$ is a $\varphi$-fibre,
so $-K_X \cdot l_i=2$. 
Hence there are two types of degenerations, either

a) $-K_X \cdot l_0=1$ and the degeneration is of the form $l_0+l_1+l_2$ with possibly $l_1=l_2$, or

b)  $-K_X \cdot l_0=3$ and the degeneration is of the form $l_0+l_1$.

{\em Case a)} Then  $-K_X \cdot l_0=1$ implies $\xi \cdot l_0=-2$, so $l_0$ is contained in $\xi$ and does not pass through $x$. Thus up to renumbering $x \in l_1$, {\em i.e.}\  
$l_1= \fibre{\varphi}{\varphi(x)}$. If we fix $l_0$, then $l_2$ varies in a 1-dimensional family parametrised by the line $\varphi(l_0)$. 
The surface covered by these degenerations is of course nothing else than $X_{\varphi(l_0)}:=\fibre{\varphi}{\varphi(l_0)}$.
One easily computes that 
$l_0 \subset X_{\varphi(l_0)}$
is a section with self-intersection $-3$, hence the surface $X_{\varphi(l_0)} \simeq \PP(E|_{\varphi(l_0)})$ is a Hirzebruch surface $\mathbb F_3$. 
A look at the extension defining $E$ shows that we are in this case exactly when $\varphi(l_0)$ meets the line $d_0$.
The restriction of the foliation $\sF_x$ to $X_{\varphi(l_0)}$ is the foliation defined by the pencil
\[
|\sO_{X_{\varphi(l_0)}}(l_0+l_1+l_2) \otimes \sI_x|.
\]
The singular locus of this foliation contains obviously the exceptional section $l_0$ but also the curve $\fibre{\varphi}{\varphi(x)}$, since the cyle
$l_0+2l_1$ is not reduced along $l_1=\fibre{\varphi}{\varphi(x)}$.

{\em Case b)} As in Case a) we show that $l_0 \subset \xi$ and $l_1= \fibre{\varphi}{\varphi(x)}$. 
If we fix $l_0$ this shows that the surface $X_{\varphi(l_0)}:=\fibre{\varphi}{\varphi(l_0)}$ contains a unique degeneration.
Moreover $l_0 \subset X_{\varphi(l_0)}$
is the exceptional section of the Hirzebruch surface 
$X_{\varphi(l_0)} \simeq \PP(E|_{\varphi(l_0)}) \simeq \mathbb F_1$ 
and all this happens exactly when $\varphi(l_0)$ does not meet the line $d_0$.
The restriction of the foliation $\sF_x$ to $X_{\varphi(l_0)}$ is the foliation defined by the pencil
\[
|\sO_{X_{\varphi(l_0)}}(l_0+l_1) \otimes \sI_x|.
\]
The general member of this pencil is a section with self-intersection one, 
the unique non-smooth member is the degeneration $l_0+l_1$.
Thus the singular locus of the foliation consists of the point $x$ and $l_0 \cap l_1=\fibre{\varphi}{\varphi(x)} \cap \xi$.
In this case the singular locus of the foliation does not depend on the choice of $l_0$.

The singular locus of $\sF_x$ is obtained as the union of the singular loci of the restricted foliations, this shows the claim.
\end{example}

\section{Enumeration of quasi-lines}

The main object of this section is to prove our only result 
when the Assumption \ref{assumption} fails 
(which corresponds to the case (II) described in Proposition \ref{propositiontwofoliations}).
As we already said, we have no concrete example but we realised
however that we can provide some information concerning
another important problem, namely bounding $e(X,l)$. 
Finally in Proposition \ref{propositionbadestimate} we try to address the same enumerative problem under the hypothesis
that Assumption \ref{assumption} holds. We will see that although Proposition \ref{propositioncorankone} provides additional
structure information, the bound obtained is far from being optimal.

Before we come to the proof, recall the following consequence of the Hodge index theorem: 
let $S$ be a projective surface and $H$ an ample divisor on $X$. If $D$ is any divisor, then
\[
(D^2) (H^2) \leq (D \cdot H)^2
\]
(see \cite[V., Ex.1.7]{Har77}).
This obviously implies: let $X$ be a projective threefold 
and $H$ be an ample divisor on $X$. If $D$ is any divisor, then
\begin{equation} \label{hodgeinequality}
(D^2 \cdot H) (H^3) \leq (D \cdot H^2)^2.
\end{equation}

\begin{proof}[Proof of Theorem \ref{maincaseII}]
Fix two general points $x$ and $y$ in $X$, and denote as usual
by $F_{x,y}$ the unique $\sF_x$-leaf through $y$
and by $F_{y,x}$ the unique $\sF_y$-leaf through $x$.

By Proposition \ref{propositiontechnicalpropertiesx},b) applied 
to the foliations $\sF_x$ and $\sF_y$,  
the quasi-lines through $x$ and $y$ are  contained in $F_{x,y}$ and $F_{y,x}$.
By hypothesis, the intersection 
$F_{x,y} \cap F_{y,x}$ is a strict subset of $F_{x,y}$ and $F_{y,x}$.
Since the rank of the foliations is two, the intersection is
a union of curves that contains the quasi-lines passing through $x$ and $y$.
It follows that 
\[
e(X,l) \deg l  \leq (F_{x,y} \cap F_{y,x}) \cdot H \leq F_{x,y} \cdot F_{y,x} \cdot H.
\]
The varieties $F_{x,y}$ 
and $F_{y,x}$ are fibres of the equidimensional fibration $q'$ (see Diagram \eqref{basicdiagramglobal}),
so they have the same cohomology class. 
By Formula \eqref{hodgeinequality} above
\[
F_{x,y}^2 H \leq \frac{(F_{x,y} H^2)^2}{H^3},
\]
so we obtain
\[
e(X,l) \leq \frac{(\deg F_{x,y})^2}{\deg X \times \deg l}.
\]
Thus we are left to bound the degree of $F_{x,y}$ : 
let \holom{\mu}{F'_{x,y}}{F_{x,y}} be the desingularisation of $F_{x,y}$ from Proposition \ref{propositiontwofoliations},b). Then we have
\[
\deg F_{x,y} = H|_{F_{x,y}}^2 = (\mu^* H)^2.
\]
Let $l_1$ be a quasi-line through $x$ and $y$. Since
$y$ is a smooth point of $F_{x,y}$, its strict transform $l_1'$ is a rational curve that passes through $y$. 
By  Proposition \ref{propositiontwofoliations},b) we know that
$y$ is in the free locus of $F'_{x,y}$, so $l_1'$ is even a free rational curve of $\mu^* H$-degree exactly $\deg l$.

By construction of the leaf $F_{x,y}$ there exists a family of quasi-lines $l_2$ through $x$ that meets $l_1$.
The strict transforms $l_2'$ of such general quasi-line form a dominant family of rational curves that meet $l_1'$ and
also have  $\mu^* H$-degree exactly $\deg l$. 
Thus by comb-smoothing we can construct a dominant 
family of rational curves of degree $2 \deg l$ passing through $y$. 
Therefore by \cite[V., Prop.2.9]{Kol96}
\[
(\mu^* H)^2 \leq 4 (\deg l)^2.
\]  
This implies the statement. 
\end{proof}

\begin{remark}
In the preceding proof, we used the hypothesis on the dimension to compute the degree of the intersection 
$F_{x,y} \cap F_{y,x}$ {\em via} 
an intersection product $F_{x,y} \cdot F_{y,x}$. If $X$ has 
higher dimension this is no longer possible,
since two surfaces meeting in a bunch of curves are a case 
of excess intersection.
\end{remark}

\begin{proposition} \label{propositionbadestimate}
Let $X$ be a projective manifold of dimension three 
containing a quasi-line $l$.
Let $x$ be a general point in $X$, and 
denote by $\sF_x$ the corresponding foliation.
Assume that $\rk \sF_x=2$ and the 
Assumption \ref{assumption} {\bf holds}.
Denote by \holom{\varphi}{\tilde X}{\PP^2} a holomorphic model of the fibration constructed
in Proposition \ref{propositioncorankone}.
Then 
\[
e(X,l) \leq 
\binom{\frac{(d+1)(d+2)}{2} \max(d, H^2)}{\frac{d(d+3)}{2}}^{\frac{d^2(d+1)(d+2)}{2}}
\]
where $d=H \cdot l$ and $H^2:=H^2 \cdot \varphi^* \sO_{\PP^2}(1)$
and $H$ is a very 
ample\footnote{Actually it is sufficient that the restriction of $H|_\sL$ to a general $\sF_x$-leaf $\sL$ is very ample.} 
line bundle on $\tilde{X}$.
\end{proposition}

\begin{proof}
Let $\sL$ be a general $\sF_x$-leaf. If $y \in \sL$ is a general point, then by Lemma \ref{lemmasecondcasequasilines},b)
any quasi-line through $x$ and $y$ is contained in the smooth locus of $\sL$ and is a quasi-line in $\sL$. 
Since $\sL$ is a surface, we have $e(\sL, l)=1$ for any quasi-line $l \subset \sL$. 
Thus bounding $e(X,l)$ is equivalent to bounding the number of families of quasi-lines on $\sL$.
By \cite[Prop. 2.1]{IV03} we have
$$
h^0(\sL, \sO_\sL(H)) \leq \frac{(d+1)(d+2)}{2}.
$$
Thus the very ample line bundle $H|_{\sL}$ defines an embedding of $\sL \hookrightarrow \PP^N$
with $N \leq \frac{d(d+3)}{2}$. Moreover the embedded surface $\sL$ has degree
$$
\deg \sL = (H|_\sL)^2 = H^2 \cdot \sL = H^2 \cdot \varphi^* \sO_{\PP^2}(1),
$$
since by Proposition \ref{propositioncorankone} the surface $\sL$ is the preimage of a line in $\PP^2$.
The statement now follows from a general estimate on the number of irreducible components of the Chow scheme
\cite[Prop.3.6]{Hei05}
\end{proof}

\begin{example} \label{examplepardinicount}
We will now bound $e(X,l)$ for the family of quasi-lines constructed in Example \ref{examplepardini}.
The proof of Proposition \ref{propositionbadestimate} shows that this is equivalent to bounding the number
of families of quasi-lines in the surface $S_d$ which is the blow-up of $\PP^2$
in five general points. The anticanonical divisor $-K_{S_d}$ is very ample, 
so Proposition \ref{propositionbadestimate} yields
$$
e(X,l) \leq \binom{60}{9}^{90}.
$$
Since we know that $h^0(S_d, -K_{S_d})=5$, the proof of the proposition shows that we can lower the bound to
$$
e(X,l) \leq \binom{20}{4}^{45}
$$
which is still huge. We will now compute explicitly the number of families on $S_d$:
if $l \subset S_d$ is a quasi-line, the linear system $|\sO_{S_d}(l)|$ is base-point free and defines a birational map $\holom{\varphi_l}{S_d}{\PP^2}$
such that $l$ maps onto a line. Since $-K_{S_d}$ is ample, we see that $\varphi_l$ 
is the blow-up of five points in general position in $\PP^2$,
in particular $\varphi_l$ contracts five disjoint $(-1)$-curves on $S_d$. 
{\em Vice versa} five disjoint $(-1)$-curves on $S_d$ determine a 
representation of $S_d$ as a blow-up of five general points in $\PP^2$.
Thus the problem reduces to counting the number of choices of five disjoint $(-1)$-curves on $S_d$.

Consider now the representation \holom{\mu}{S_d}{\PP^2} fixed in Example \ref{examplepardini}, and denote by
$p_1, \ldots, p_5 \in \PP^2$ the points we blow up. With respect to this representation, the $(-1)$-curves
on $S_d$ are 
\begin{itemize}
\item the five exceptional divisor $E_1, \ldots, E_5$,
\item the (strict transforms of the) ten lines passing through exactly two of the points $p_1, \ldots, p_5$
which we denote by $d_1, \ldots, d_{10} \subset S_d$, and
\item the (strict transform of the) unique conic passing through $p_1, \ldots, p_5$ which we denote by $C \subset S_d$.
\end{itemize}
Out of these sixteen curves we have to choose a configuration of five that are disjoint. 

{\it Case 1) The configuration contains $C$.} The conic meets all the $E_i$ and is disjoint from all the $d_j$. 
Thus we are left to choose a configuration of four disjoint $d_j$. The dual graph of the $d_j$ is the famous Petersen graph
and it is easy to see that there are exactly five such possibilities.
\begin{figure}[h]
\begin{center}
\includegraphics[height=30mm]{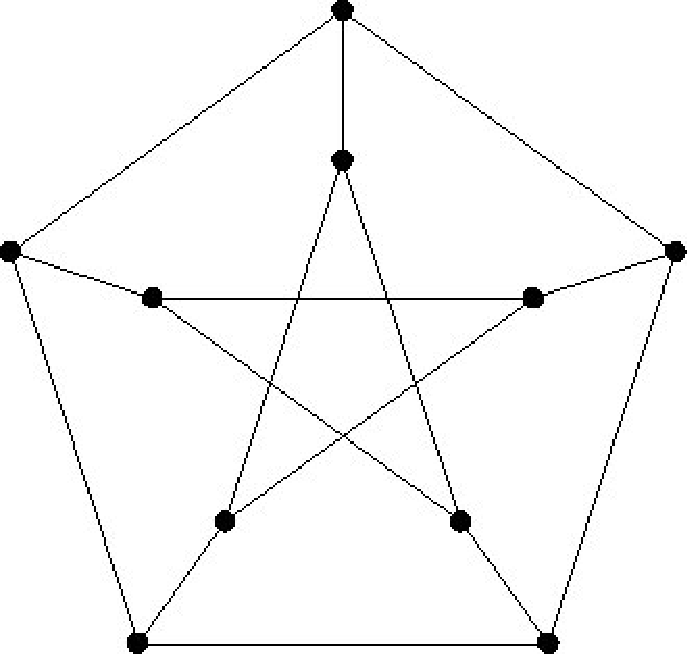}
\end{center}
\end{figure}

{\it Case 2) The configuration does not contain $C$.}
The Petersen graph shows that it is not possible to choose five disjoint $d_j$.
If the configuration contains an exceptional divisor $E_i$, this excludes the four $d_j$ coming from lines through $p_i$. 
A look at the Petersen graph shows that it is not possible to choose four disjoint $d_j$ among the remaining six,
so the configuration contains at least two exceptional divisors $E_i, E_{i'}$ which excludes seven $d_j$ coming from lines 
through $p_i$ or $p_{i'}$. The remaining three $d_j$ are disjoint, so we get ten configurations
with exactly two exceptional divisors. In the same way we see that if there are at least three exceptional divisors
in the configuration, the configuration is given by $E_1, \ldots, E_5$.

In total we get $5+10+1$ configurations, so
$$
e(X,l) \leq 16.
$$
\end{example}

\end{document}